\documentclass[12pt]{article}
\title{LEBESGUE SPACES AND INCLUSION RELATION}
\author{{\sc C. Ganesamoorthy}\\Department of mathematics \\Alagappa University\\Karaikudi - 630003\\India.\\Email: ganesamoorthyc@gmail.com, gmoorthyc@alagappauniversity.ac.in}
\date{}
\usepackage{amsmath, wasysym, mathrsfs, amssymb, yfonts}
\begin{document}
\maketitle
\newtheorem{lem}{\sf Lemma}[section]
\newtheorem{defn}[lem]{\sf Definition}
\newtheorem{thm}[lem]{\sf Theorem}
\newtheorem{cor}[lem]{\sf Corollory}
\newtheorem{pro}[lem]{\sf Property}
\newtheorem{exm}[lem]{\sf Example}
\begin{abstract}
The proofs of A. Villani on inclusion relations among classical Lebesgue spaces are dicussed. The techinque of using closed graph theorem, due to Villani, is applied to derive results on inclusion relations among some more additional spaces. An application of automatic continuity of multiplicative linear functionals on Banach algebras is found.
\end{abstract}
\section{Introduction}

\hspace{.3cm} Let $(X, \mathscr{M}, \mu)$ be a positive measure space. For $1\leq p \leq \infty$, let $(L^p, \|.\|_p)$ denote the classical Lebesgue space derived from $(X, \mathscr{M}, \mu)$. Each $(L^p, \|.\|_p)$  is a Banach space. The natural question about the inclusion relation is: when does the relation $L^p \subset L^q$ hold, for $p<q $ or $q<p$ ? This question appears in some text books. This question also appears in some research notes (see: \cite{ro}, \cite{bs} and \cite{vil}). The answer of this question appears in such notes, because of the techniques applied to derive the answer. One note of A. Villani \cite{vil} applies the closed graph theorem to derive the answer. So, the first part of the article is to be present the work of A. Villani \cite{vil}, which uses the closed graph theorem of functional analysis. The second part of this article presents an application of the following result to obtain an extended answer: Every multiplicative linear functional on a Banach algebra is continuous. Let $\mathscr{M}_0=\{E\in \mathscr{M}:\mu(E)>0\}$ and $\mathscr{M}_{\infty}=\{E\in \mathscr{M}:\mu(E)<\infty\}$. The following two theorems are known: \\
\begin{thm} \label{thm0.1}
The following are equivalent: \begin{description}
\item(i) $L^p \subset L^q $ for some $p, q \in [1, \infty]$ with $p<q$
\item(ii) $\inf\limits_{E \in \mathscr{M}_0} \mu(E)>0$
\item(iii) $L^p \subset L^q$ for all $p, q \in [1, \infty]$ with $p<q$
\end{description}
\end{thm}

\begin{thm} \label{thm0.2}
The following are equivalent: \begin{description}
\item(i) $L^p \supset L^q $ for some $p, q \in [1, \infty)$ with $p<q$
\item(ii) $\sup\limits_{E \in \mathscr{M}_{\infty}} \mu(E)<\infty$
\item(iii) $L^p \supset L^q$ for all $p, q \in [1, \infty)$ with $p<q$
\end{description}

\end{thm}

The next section presents proofs of Villani \cite{vil} for these theorems based  on the closed graph theorem. Observe that the case of assuming $\infty$ is not included in the Theorem \ref{thm0.2}. We shall need the following two results (see: \cite{rc}  \cite{rjl}, \cite{ij}, \cite{rl}, \cite{ol}, \cite{rj}, \cite{yr} and \cite{ew} ) to apply the closed graph theorem.

\begin{lem} 
\label{lem0.3} If $1\leq p\leq \infty$ and if $(f_n)$ is a Cauchy sequence in $L^p$ with limit $f$, then $(f_n)$ has a subsequence which converges pointwisely to $f$ almost everywhere on $X$.
\end{lem}
\begin{lem} \label{lem0.4}
If $(f_n)$ is a cauchy sequence in $L^{\infty}$ with limit $f$, then $(f_n)$ converges pointwisely to $f$ almost everywhere on $X$.
\end{lem}

\section{Standard inclusion relations}

\hspace{.3cm} This section presents proofs of Villani \cite{vil} for the Theorems \ref{thm0.1}, and \ref{thm0.2}, which are based on the closed graph theorem.
\begin{lem} \label{lem1.1}
Let $p, q \in [1, \infty]$. The set theoretic inclusion $L^p \subset L^q$ implies that the inclusion map $I:L^p \longrightarrow L^q$ is continuous. \end{lem} 
{\sf Proof:} Let $(f_n)$ be a sequence which converges to $f$ in $L^p$; and suppose $(f_n) $ converges to $g$ in $L^q$. Then, by Lemma \ref{lem0.3}, $(f_n)$ has a subsequence $(g_m)$ which converges pointwisely to $f$ almost everywhere. The same reason implies that $(g_m) $  has a subsequence $(h_k)$ which converges pointwisely to $g$ almost everywhere. So, $(h_k)$ converges pointwisely to $f$ as well as to $g$ almost everywhere. So $f=g$ almost everywhere. This shows that the inclusion map $I$ has closed graph. So, the closed graph theorem implies the continuity of $I$.\hfill $\blacksquare$\\

So, we applied the closed graph theorem in our development.\\

\noindent{\sf Proof of the Theorem \ref{thm0.1}:}  $(i)\Rightarrow(ii)$: Suppose that (i) is true. By Lemma \ref{lem1.1}, the inclusion map $I:L^p \longrightarrow L^q$ is continuous. So, there is a positive constant $k$ such that $\|f\|_q \leq k\|f\|_p$, for all $f \in L^p$. For every $E \in \mathscr{M}_0 \cap \mathscr{M}_{\infty}$, the characteristic function $\chi_E \in L^p \cap L^q$, and hence $\|\chi_E\|_q \leq k\|\chi_E\|_p$.\\ Case (a): Suppose $q<\infty$. Then, for every $E\in \mathscr{M}_0 \cap \mathscr{M}_\infty$, we have $\mu(E)^{1/q} \leq k \mu(E)^{1/p}$ and hence $0< $ {\Large $\left(\frac{1}{k}\right)^{pq/q-p}$} $\leq \mu(E)$. Therefore, $\inf\limits_{E\in \mathscr{M}_0} \mu(E)>0$.\\ Case (b): Suppose $q=\infty$. Then for every $E\in \mathscr{M}_0 \cap \mathscr{M}_\infty$, we have $1\leq k\mu(E)^{1/p}$ and hence {\Large $\left(\frac{1}{k}\right)^{p}$} $\leq \mu(E)$. So, again we have $\inf\limits_{E\in \mathscr{M}_0} \mu(E)>0$.\\
$(ii)\Rightarrow (iii)$: Suppose $(ii)$ is true. Let $1\leq p<\infty$. Let $f\in L^p$. Let $E_n =\{x\in X: |f(x)|>n\}$, for $n=1,2, \cdots$ . Then $\int\limits_{X} |f|^p \chi_{E_n} d\mu \geq \int\limits_{X} n^p \chi_{E_n} d\mu$ and hence $\int\limits_{E_n} |f|^p d\mu \geq  n^p \mu(E_n) $ for every $n$. Thus $0\leq\mu(E_n)\leq$ {\Large $\frac{\int\limits_{E_n} |f|^p d\mu}{n^p}$} $\leq$ {\Large $\frac{\int\limits_{X} |f|^p d\mu}{n^p}$} $\rightarrow 0$ as $n\rightarrow \infty$. So, $\mu(E_n)\rightarrow 0$ as $n\rightarrow \infty$. By condition $(ii)$, there is an integer $n_0$ such that $\mu(E_n)=0$ for every $n\geq n_0$. Therefore $|f|\leq n_0$ almost everywhere on $X$. So,  $\|f\|_{\infty} \leq n_0 <\infty$, and hence $f\in L^{\infty}$. That is, $(iii)$ is true whenever $q=\infty$. Suppose $p<q<\infty$. Let $E=\{x\in X:|f(x)|\leq 1\}$ and $F=\{x\in X:|f(x)|> 1\}$. Since $\mu(F) \leq \int\limits_{F} |f|^p d\mu < \infty$, we have $\int\limits_{F} |f|^q d\mu \leq \|f\|_{\infty}^q \mu(F) <\infty$. So, $\int\limits_X |f|^q d\mu = \int\limits_E |f|^q d\mu + \int\limits_F |f|^q d\mu \leq \int\limits_E |f|^p d\mu + \int\limits_F |f|^q d\mu < \infty$. Thus, if $f\in L^p$, then $f\in L^q$, whenever $p<q$. This proves $(iii)$ for all cases.\\
$(iii) \Rightarrow (i)$: Obvious.\hfill $\blacksquare$\\

\noindent {\sf Proof of Theorem \ref{thm0.2}:}
$(i)\Rightarrow (ii)$: Suppose that $(i)$ is true. By Lemma \ref{lem1.1}, the inclusion map $I:L^q \rightarrow L^p$ is continuous. So, there is a positive constant $k$ such that $\|\chi_E\|_p \leq k\|\chi_E\|_q$, for every $E\in \mathscr{M}_0 \cap \mathscr{M}_\infty$. So, for every $E\in \mathscr{M}_0 \cap \mathscr{M}_\infty$, we have $\mu(E)^{1/p} \leq k\mu(E)^{1/q}$ and hence $\mu(E) \leq k^{pq/q-p} < \infty$. This proves $(ii)$. \\ $(ii) \Rightarrow (iii)$: Suppose $(ii)$ is true. Suppose $1\leq p<q<\infty$; and let $f\in L^q$. To each $n=1,2,\cdots$, define $E_n=${\Large$\{$}$x\in X$:{\Large $\frac{1}{n+1}$}$\leq |f(x)|<${\Large $\frac{1}{n}\}$}. Then \\
{\Large $\frac{1}{(n+1)^q}$}$\mu(E_n)=\int\limits_X${\Large $\frac{1}{(n+1)^q}$}
$\chi_{E_n} d\mu \leq \int\limits_X |f|^q \chi_{E_n} d\mu \leq \int\limits_X |f|^q d\mu < \infty$. In particular \\$\mu(E_n)<\infty$, for every $n$. By condition $(ii)$, we have $\sum\limits_{n=1}^{\infty} \mu(E_n) = \sup\limits_n \mu(E_1 \cup E_2 \cup \cdots\cup E_n) < \infty$. Define $E=\{x\in X: 0< |f(x)| < 1\}$ and $F=\{x\in X:|f(x)| \geq 1\}$. Then $E=\mathop{\cup}\limits_{n=1}^{\infty} E_n$. So, $\mu(E)=\mu(\mathop{\cup}\limits_{n=1}^{\infty} E_n)=\sum\limits_{n=1}^{\infty} \mu(E_n)<\infty$. This implies that $\int\limits_X |f|^p d\mu = \int\limits_E |f|^p d\mu + \int\limits_F |f|^p d\mu \leq \int\limits_E |f|^p d\mu +\int\limits_F |f|^q d\mu < \infty$. Thus $f\in L^p$ and hence $L^q \subset L^p$.\\ $(iii)\Rightarrow (i)$: Obvious. \hfill $\blacksquare$\\

\section{$L^{\infty}$ case}
\begin{exm} \label{exm2.1}
Let $X=\{a, b\}; \ \mathscr{M} = \{\phi,  \{a\},  \{b\},  X\}$. Define $\mu$ by $\mu(\phi)=0, \ \mu(X) = \mu(\{b\})= \infty$ and $\mu(\{a\})=1$. Then $\sup\limits_{E\in \mathscr{M}_{\infty}} \mu(E)=1 <\infty$. Define $f:X\rightarrow \{0, 1\}$ by $f(a)=0$ and $f(b) =1$. Then $f\notin L^1$ but $f\in L^{\infty}$. So, $L^{\infty}$ is not contained in $L^1$. 
\end{exm} 

This shows that the Theorem \ref{thm0.2} cannot be extended to the case $q=\infty$. However, we have the following theorem
\begin{thm} \label{thm2.2}
The following are equivalent:
\begin{description}
\item(i) $L^p \supset L^{\infty} $ for some $p\in [1, \infty)$ 
\item(ii) $\mu(X)<\infty$
\item(iii) $L^p \supset L^{\infty}$ for all $p \in [1, \infty)$ 
\end{description}
\end{thm}
{\sf Proof:} $(i)\Rightarrow (ii)$: Suppose $(i)$ is true. By Lemma \ref{lem1.1}, there is a  constant $k>0$ such that $\|f\|_p \leq k\|f\|_{\infty}$ for every $f\in L^{\infty}$. In particular, we have $\|\chi_X\|_p \leq k\|\chi_X\|_{\infty}$. Thus $\mu(X) \leq k^p <\infty$.\\$(ii)\Rightarrow (iii)$: Suppose $(ii)$ is true. Let $1\leq p<\infty$ and let $f\in L^\infty$. Then $\int\limits_X |f|^p d\mu \leq \|f\|_{\infty}^p \int\limits_X d\mu = \|f\|_{\infty}^p \mu(X)<\infty$. Thus $f\in L^p$. This shows that $L^p \supset L^{\infty}$. \\$(iii)\Rightarrow(i)$: Obvious.\hfill $\blacksquare$\\

Can we replace $L^{\infty}$ by some other spaces so that Theorem \ref{thm0.2} can be extended to those spaces ? The answer is `yes'. One such space is the dual of $L^1$ (See: \cite{sc}). However, another simple space will be discussed here.\\ 

Each $L^p, 1\leq p\leq \infty$, is a Banach space in which elements are equivalence classes. The equivalence relation is: $f$ is equivalent to $g$ if and only if $f=g$ almost everywhere on $X$. Let us consider another equivalence relation `$\AC$' defined by the rule: $f\AC g$ if and only if $f=g$ almost everywhere on $E$, for every $E\in \mathscr{M}_{\infty}$. This equivalence relation $\AC$ defines new equivalence classes. However, they are same (identified) as the usual classes in $L^p$, for $1\leq p< \infty$. But, they may be different equivalence classes in $L^\infty$ (refer to Example \ref{exm2.1}). When the elements of $L^\infty$ are identified by the new equivalence classes induced by `$\AC$', let us write $L^{\infty}_{qt}$ to denote the corresponding space with the norm $\|.\|_{qt\infty}$ defined by $\|f\|_{qt\infty} = \sup\limits_{E\in \mathscr{M}_{\infty}}  \|f\chi_E\|_{\infty}$. Here $f$ is identified with the equivalence class in which $f$ is a member. Then $(L_{qt}^{\infty}, \|.\|_{qt\infty})$ is a Banach space in view of the following reasons: If $N=\{f\in L^\infty : f\chi_E =0$ almost everywhere on $E$ for every subset $E\in \mathscr{M}_{\infty}\}$, then $N$ is a closed subspace of $(L^{\infty}, \|.\|_{\infty})$; in the view of Lemma \ref{lem0.4}. Then $(L^{\infty}_{qt}, \|.\|_{qt\infty})$ may be realized as the quotient space $L^{\infty}/N$ with the quotient norm.\\

One can derive the following Lemma \ref{lem2.3} which is similar to the Lemma \ref{lem1.1}

\begin{lem} \label{lem2.3}
Let $p\in [1, \infty)$.
\begin{description}
\item (a) The set theoretic inclusion $L^p \subset L^{\infty}_{qt}$ implies that the inclusion map $I:L^p \longrightarrow L^{\infty}_{qt}$ is continuous. 
\item (b) The set theoretic inclusion $L^{\infty}_{qt} \subset L^p$ implies that the inclusion map $I:L^{\infty}_{qt} \longrightarrow L^p$ is continuous. 
\end{description}
\end{lem}  

Now we shall use this Lemma \ref{lem2.3} to derive the following theorem.
\begin{thm} \label{thm2.4}
The following are equivalent:
\begin{description}
\item(i) $L^p \supset L^{\infty}_{qt} $ for some $p\in [1, \infty)$ 
\item(ii) $\sup\limits_{E\in \mathscr{M}_{\infty}} \mu(E)<\infty$
\item(iii) $L^p \supset L^{\infty}_{qt}$ for all $p \in [1, \infty)$ 
\end{description}
\end{thm}
{\sf Proof:} $(i)\Rightarrow (ii)$: Suppose $(i)$ is true. Then the inclusion map $I:L^{\infty}_{qt} \longrightarrow L^p$ is continuous. So, there is a  constant $k>0$ such that $\|f\|_p \leq k\|f\|_{qt\infty}$, for every $f\in L^{\infty}_{qt}$. In particular, for every $E \in \mathscr{M}_0\cap \mathscr{M}_{\infty}$, we have $\|\chi_E\|_p \leq k\|\chi_E\|_{qt\infty}; \ \mu(E)^{1/p}\leq k$ and hence $\mu(E) \leq k^p <\infty$.\\$(ii)\Rightarrow (iii)$: Suppose $(ii)$ is true. Let $1\leq p < \infty$, and let $f \in L^{\infty}_{qt}$. Then there is an increasing sequence $(E_n)$ in $\mathscr{M}_{\infty}$ such that $\|f\|_{qt\infty}= \sup\limits_n \|f\chi_{E_n}\|_{\infty}$ and such that $\mu((X\setminus \mathop{\cup}\limits_{n=1}^\infty E_n)\cap F)=0$ for every $F\in \mathscr{M}_\infty$. Write $E=\mathop{\cup}\limits_{n=1}^{\infty} E_n$. Define $g$ on $X$ by $g(x) = \begin{cases} f(x) &\text{ if } x\in E\\0 & \text{ if } x \notin E \end{cases}$ \\ Then $f\AC g$. Also $\mu(E)=\sup\limits_n \mu(E_n)< \infty$. So,  $\int\limits_X |g|^p d\mu \leq \|f\|_{qt\infty}^p \mu(E)<\infty$. This shows that $L^p \supset L^{\infty}_{qt}$. \\$(iii)\Rightarrow(i)$: Obvious.\hfill $\blacksquare$\\

$L_{qt}^{\infty} \supset L^p$ if and only if $L^{\infty}\supset L^p$, for any fixed $p\in [1, \infty)$. So, Theorem \ref{thm0.1} gives the following theorem.
\begin{thm} \label{thm2.5}
The following are equivalent: \begin{description}
\item(i) $L^p \subset L^{\infty}_{qt} $ for some $p\in [1, \infty)$ 
\item(ii) $\inf\limits_{E \in \mathscr{M}_0} \mu(E)>0$
\item(iii) $L^p \subset L^{\infty}_{qt}$ for all $p \in [1, \infty)$
\end{description}
\end{thm}

One can also give an independent proof based on the arugments used to establish the Theorem \ref{thm0.1}. So, we have spaces to fill up the gap (case: $\infty$) in the Theorem \ref{thm0.2}. But, we achieved this through a new equivalence relation `$\AC$' defined above. The usual equivalence relation can also be used to fill up the gap (case: $\infty$) in the Theorem \ref{thm0.2}. This is achieved in the next secion. In the next section two functions are equivalent if and only if they are equal almost everywhere on the entire set $X$.  

\section{The space $C_0$}

\begin{defn}
$C_0 =\{f\in L^{\infty}$: For every $\epsilon >0$, there is a set $E\in \mathscr{M}_{\infty} $ such that $ \|f\chi_{X\setminus E}\|_{\infty}<\epsilon\}$. 
\end{defn} 

This $C_0$ becomes a closed subspace of $(L^{\infty}, \|.\|_{\infty})$, when members of $C_0$ are usual equivalence classes. So, $(C_0, \|.\|_{\infty})$ is a Banach space.

\begin{thm}
The following are equivalent:
\begin{description}
\item(i) $L^p \supset C_0 $ for some $p\in [1, \infty)$ 
\item(ii) $\sup\limits_{E\in \mathscr{M}_{\infty}} \mu(E)<\infty$
\item(iii) $L^p \supset C_0$ for all $p \in [1, \infty)$ 
\end{description}
\end{thm}
{\sf Proof:} $(i)\Rightarrow (ii)$: The arugment used for Theorem \ref{thm2.4} ($(i)\Rightarrow (ii)$)is applicable. \\ $(ii)\Rightarrow (iii)$: Since $C_0$ can be identified as a subspace of $L^{\infty}_{qt}$, Theorem \ref{thm2.4} can be applied to derive this implication.\\ $(iii)\Rightarrow(i)$: Obvious.\hfill $\blacksquare$\\

\begin{thm} \label{thm3.3}
The following are equivalent: \begin{description}
\item(i) $L^p \subset C_0$ for some $p\in [1, \infty)$ 
\item(ii) $\inf\limits_{E \in \mathscr{M}_0} \mu(E)>0$
\item(iii) $L^p \subset C_0$ for all $p \in [1, \infty)$
\end{description}
\end{thm}
{\sf Proof:} $(i)\Rightarrow(ii)$: Since $C_0\subset L_{qt}^{\infty}$; Theorem \ref{thm2.5} proves this implication.\\$(ii)\Rightarrow(iii)$: Suppose $(ii)$ is true. Let $1\leq p<\infty$. Let $f \in L^p$. Theorem \ref{thm0.1} implies that $f \in L^{\infty}$. Let $E_n=\{x\in X:|f(x)|>${\Large $\frac{1}{n}$}$\}$ , for $n=1,2\cdots$. Then $\int\limits_{E_n} |f|^p d\mu \geq ${\Large$\frac{1}{n^p}$}$\mu(E_n)$ and hence $E_n \in \mathscr{M}_{\infty}$, for every $n$. Also, $\|f\chi_{X\setminus E_n}\|_{\infty} \rightarrow 0$ as $n\rightarrow \infty$. So, $f\in C_0$, and hence $(iii)$ is true.\\$(iii)\Rightarrow (i)$: Obvious \hfill $\blacksquare$

\section{Dual and Predual of $L^1$}

\hspace{.3cm} The conditions $(ii)$ of earlier theorems can be changed into other forms. Two non trivial forms are presented in this section. We shall use the result: Every multiplicative linear functional on a Banach algebra is continuous.\\

\begin{thm}
The following are equivalent: 
\begin{description}
\item(i) Dual of $C_0$ is (isomorphic with) $L^1$ 
\item(ii) $\inf\limits_{E\in \mathscr{M}_0} \mu(E) >0$.
\end{description}
\end{thm}
{\sf Proof:} $(ii) \Rightarrow(i)$: Suppose $(ii)$ is true. Let $T$ be a non zero continuous linear functional on $C_0$. For each atom $E\in \mathscr{M}_0 \cap \mathscr{M}_{\infty}$, let $a_E = T(\chi_E)$. The collection of all atoms $E\in \mathscr{M}_0 \cap \mathscr{M}_{\infty}$ with $a_E\neq 0$ should be countable; and $\sum |a_E|<\infty$ where the sum runs over all $E$ in this collection. Otherwise, there is a sequence $E_1, E_2, \cdots $ of mutually disjoint atoms in  $\mathscr{M}_0 \cap \mathscr{M}_{\infty}$ such that $\sum\limits_{i=1}^{\infty} |a_{E_i}|=\infty$ and $a_{E_i} \neq 0$ for every $i$; and in this case the sequence of functions \\$f_n(x) = \begin{cases}$ {\Large $\frac{\bar{a}_{E_i}}{|a_{E_i}|}$}$ & \text{ if } x\in E_i \text{ for some } i=1,2,\cdots,n \\ \ \ 0 &\text{ if } x\notin E_i \text{ for all } i=1,2,\cdots,n \end{cases}$ \\is bounded in $C_0$ and the sequence $T(f_n)=\sum\limits_{i=1}^{n} |a_{E_i}|$ is not bounded. 

Let $F_1, F_2, \cdots$ or $F_1, F_2, \cdots, F_k$ be the set of all atoms in $\mathscr{M}_0 \cap \mathscr{M}_{\infty}$ such that $a_{F_i} \neq 0$. We lose no generality by assuming that these atoms are mutually disjoint. Define $f$ on $X$ by $f(x) = \begin{cases}$ {\Large $\frac{a_{F_i}}{\mu(F_i)}$}$ & \text{ if } x\in F_i \text{ for some } i\\ \ \ 0 &\text{ if } x\notin F_i \text{ for all } i\end{cases}$ .\\ Then $\int\limits_X |f| d\mu = \sum_{i=1}^{\infty \text{ (or) } k} |a_{F_i}|<\infty$. Thus $f\in L^1$. If $g\in C_0$, then, to each $i$, there is a number $b_i$ such that $g$ assumes $b_i$ almost everywhere in $F_i$; and hence $\int\limits_X fg d\mu = \sum_{i=1}^{\infty \text{ (or) } k} a_{F_i}b_i = T(g) $. This is an application of $(ii)$ and continuity of $T$. It is easy to verify that $\|T\|=\|f\|_1$. 

On the other hand, every $f\in L^1$ corresponds to a continuous linear functional $T$ on $C_0$ such that $T(g) = \int\limits_X  fg d\mu$, for all $g \in C_0$; and it is easy to verify that $\|T\|=\|f\|_1$.  \\
$(i)\Rightarrow (ii)$: Suppose $(ii)$ is not true. Then there is an infinite sequence of mutually disjoint sets $E_1, E_2, \cdots \in \mathscr{M}_0 \cap \mathscr{M}_{\infty}$ such that the set $E = \mathop{\cup}\limits_{n=1}^{\infty} E_n \in \mathscr{M}_{\infty}$ and $\mu(\mathop{\cup}\limits_{i=1}^{n} E_i) \neq \mu(E)$ for every $n$. Let \textswab{U} be the family of all collections $(E_\lambda)_{\lambda \in \Lambda}$ with the properties: the index set $\Lambda$ is infinite, $E_{\lambda}\in \mathscr{M}_0 \cap \mathscr{M}_{\infty}$ for every $\lambda \in \Lambda$, $E_{\lambda} \subset E$ for every $\lambda \in \Lambda$, and no finite union of $E_{\lambda}$'s has measure $\mu(E)$. By Zorn's lemma, there is a maximal element, say, $(E_{\lambda})_{\lambda \in \Lambda}$ in \textswab{U}. Let $J'$ be the smallest ideal in $C_0$ which contains all $\chi_{E_\lambda}, \lambda \in \Lambda$. Then $J'$ is a proper ideal in the closed subalgebra $\{ f\chi_E: f\in C_0\}=A$ (say) of the Banach algebra $C_0$. Let $J''$ be a maximal ideal which contains $J'$ in the Banach algebra $A$ which contains the identity element $\chi_E$. Then this ideal $J''$ must be the kernel of a non zero continuous linear functional $T'$ on $A$ which preserves multiplication (See: [3])     

Let $T$ be the non zero continuous linear functional on $C_0$ defined by $T(f) = T'(f\chi_E)$. Suppose there is a function $f$ in $L^1$ such that $T(g)=\int\limits_X fg d\mu$ for every $g$ in $C_0$. Then we must have $T(\chi_F) = \int\limits_F fd\mu$ for every $ F \in \mathscr{M}_0 \cap \mathscr{M}_{\infty}$. Since $\{E_\lambda: \lambda\in \Lambda\}$ is maximal:
\begin{description}
\item(a) Whenever $F', F \in \mathscr{M}_0 \cap \mathscr{M}_{\infty}, \ F'\subset F$ and $F\in \{E_\lambda: \lambda\in \Lambda\}$, we have $F'\in \{E_\lambda: \lambda\in \Lambda\}$ 
\item(b) For every $F \in \mathscr{M}_0 \cap \mathscr{M}_{\infty}$ with $E\setminus F \in \mathscr{M}_0 \cap \mathscr{M}_{\infty}$ and $ F\subset E$ either $F \in  \{E_\lambda: \lambda\in \Lambda\}$ or $E\setminus F \in \{E_\lambda: \lambda\in \Lambda\}$.  
\end{description}
This shows that, for every $F \in \mathscr{M}_0 \cap \mathscr{M}_{\infty}$ with $E\setminus F \in \mathscr{M}_0 \cap \mathscr{M}_{\infty}$ and $ F\subset E$ either $\int\limits_F f d\mu =0$ or $\int\limits_{E\setminus F} fd\mu =0$, because $T(\chi_{E_\lambda}) =0$ for every $\lambda$. So, $f$ is zero almost every where on $F$ or on $E\setminus F$, for every $F\subset E$ such that $F, E\setminus F \in \mathscr{M}_0 \cap \mathscr{M}_{\infty}$ (consider positive and negative parts of real and imaginary parts of $f$) So, $f$ should be zero almost everywhere on $E$. This shows that $T$ should be zero functional; and this is a contradiction. Therefore $L^1$ cannot be the dual of $C_0$. This completes the proof. \hfill $\blacksquare$\\

\noindent We observe that the following are equivalent:
\begin{description}
\item$(i)$ There are $E_1, \ E_2 \in \mathscr{M}$ such that $X=E_1\cup E_2, \ E_1\cap E_2=\phi, \  E_1 \in \mathscr{M}_{\infty}$, and $\mu(E_2) =0$ or $E_2$ is an atom of infinite measure.
\item$(ii)$ $\sup\limits_{E\in \mathscr{M}_{\infty}} \mu(E) < \infty $
\end{description}

If this condition $(i)$ is true, then the dual of $L^1$ is $L^{\infty}_{qt}$, because $L^{\infty}_{qt}$ is equivalent to $L^{\infty}$ on $E_1$, the set mentioned in $(i)$. So, we have the following theorem.

\begin{thm}
The following are equivalent:
\begin{description}
\item(i) Dual of $L^1$ is (isomorphic with) $L^{\infty}_{qt}$ 
\item(ii) $\sup\limits_{E\in \mathscr{M}_{\infty}} \mu(E) <\infty$.
\end{description}
\end{thm}

\end{document}